\documentclass[12pt, reqno]{amsart}

\usepackage[colorlinks=true]{hyperref}
\hypersetup{urlcolor=blue, citecolor=red}

\newtheorem{theorem}{Theorem}[section]
\newtheorem{lemma}[theorem]{Lemma}

\newtheorem{corollary}[theorem]{Corollary}
\theoremstyle{definition}
\newtheorem{definition}[theorem]{Definition}

\newtheorem{remark}[theorem]{Remark}
\numberwithin{equation}{section}

\title[Fractional non-autonomous evolution equations] 
      {Cauchy problem for fractional non-autonomous evolution equations}

\author[Pengyu Chen, Xuping Zhang and Yongxiang Li]{}


 \email{chpengyu123@163.com}
 \email{lanyu9986@126.com}
\email{liyxnwnu@163.com}

\thanks{$^{\ast}$ Corresponding author.}

\begin{document}
\maketitle

\centerline{\scshape Pengyu Chen$^{\ast}$, Xuping Zhang and Yongxiang Li}
\medskip
{\footnotesize
\centerline{ Department of Mathematics, Northwest Normal University,
}
   \centerline{Lanzhou 730070, P.R. China}
}

\medskip
\bigskip

\begin{abstract}
This paper deals with the following Cauchy problem to nonlinear time fractional non-autonomous integro-differential evolution equation of mixed type via measure of noncompactness
$$
\left\{\begin{array}{ll}
  ^CD^{\alpha}_tu(t)+A(t)u(t)= f(t,u(t),(Tu)(t), (Su)(t)),\quad  t\in [0,a], \\[12pt]
   u(0)=A^{-1}(0)u_0
 \end{array} \right.
$$
in infinite-dimensional Banach space $E$, where $ ^CD^{\alpha}_t$ is the standard Caputo's fractional time derivative of order $0<\alpha\leq 1$, $A(t)$ is a family of closed linear operators defined on a dense domain $D(A)$ in Banach space $E$ into $E$ such
that $D(A)$ is independent of $t$,   $a>0$ is a constant, $f:[0,a]\times E\times E\times E\rightarrow E$ is a Carath\'{e}odory type function, $u_0\in E$,
$T$ and $S$ are Volterra and Fredholm integral operators, respectively. Combining the theory of fractional calculus and evolution families, the fixed point theorem with respect to convex-power
condensing operator and a new estimation technique of the measure
of noncompactness, we obtained the existence of mild solutions
under the situation that the nonlinear function satisfy some appropriate local growth condition and a noncompactness
measure condition. Our results generalize and improve some previous results on this topic, since the condition of uniformly continuity of the nonlinearity is not required, and also  the strong
restriction on the constants in the condition of noncompactness measure is completely deleted. As samples of applications, we consider the initial value problem to a class of time fractional non-autonomous partial differential equation with
homogeneous Dirichlet boundary condition at the end of this paper.
\end{abstract} \maketitle

\section{Introduction and main results}
\vskip5mm

Let $E$ be a real Banach space with norm $\|\cdot\|$, and let $D$ be a bounded and convex set in $E$. If an operator $\mathcal {A}:D\rightarrow E$ is completely continuous, then $\mathcal {A}$ has at least fixed point in $D$. This is the well known Schauder's fixed point theorem, which is a famous and important fixed point theorem and it has extremely widespread application. However, Schauder's fixed point theorem require the operator is completely continuous, which is a very strong restriction condition. It is well known that the solution operator of differential equations in infinite dimensional spaces are not compact, it means that Schauder's fixed point theorem can not be applied to differential equations in infinite dimensional Banach spaces. In order to overcome this strong restriction, the condition of completely continuous operator in Schauder's fixed point theorem has been relaxed to condensing operator based on the concept of condensing operator. In this way, the famous  Sadovskii's fixed point theorem has been obtained.
\vskip 3mm

Sadovskii's fixed point theorem is an important tool to study various differential equations and integral equations on infinite dimensional Banach spaces. In 1981, Lakshmikantham and Leela \cite{ll81} studied the following initial value
problem (IVP) of ordinary differential equation in Banach space $E$
\begin{equation}\label{01}
\left\{\begin{array}{ll}
u'(t)= f(t,u(t)),\quad t\in [0,a],\\[12pt]
   u(0)=u_0,
  \end{array} \right.
\end{equation}
where $a>0$ is a constant. The authors proved that, if for any constant $R>0$, $f$ is uniformly continuous
on $[0,a]\times B_R$ and satisfies the noncompactness measure condition
\begin{equation}\label{02}
\mu(f(t,U))\leq l\mu(U),\qquad \forall\;t\in [0,a],\qquad U\subset B_R,
\end{equation}\vskip2mm
\noindent where $B_R=\{u\in E:
\| u\|\leq R\}$, $l$ is a positive constant, then IVP \eqref{01} has a global solution provided that $l$ satisfies the condition
\begin{equation}\label{03}
al<1.
\end{equation}\vskip3mm

In 1989, Guo \cite{g89} studied the global solutions of the following initial
value problem (IVP) for first-order nonlinear integro-differential equations of mixed type
in a real Banach space $E$
\begin{equation}\label{04}
\left\{\begin{array}{ll}
u'(t)= f(t,u(t),(Tu)(t), (Su)(t)),\quad t\in [0,a],\\[12pt]
   u(0)=u_0,
  \end{array} \right.
\end{equation}
where
\begin{equation}\label{12}
(Tu)(t)=\int_0^tK(t,s)u(s)ds,
\end{equation}
is a Volterra integral operator with integral kernel $K\in
C(\triangle,\mathbb{R})$, $\triangle=\{(t,s)\mid 0\leq s\leq t\leq
a\}$ and
\begin{equation}\label{13}
(Su)(t)=\int_0^aH(t,s)u(s)ds,
\end{equation}
is a Fredholm integral operator with integral kernel $H\in C(\triangle_0,\mathbb{R})$, $\triangle_0=\{(t,s)\mid0\leq t, s\leq a\}$. Denote $K_0=\max\limits_{(t,s)\in \triangle}|K(t,s)|$ and $H_0=\max\limits_{(t,s)\in \triangle_0}|H(t,s)|$. Then Guo proved that IVP \eqref{04} exists at least one global solution if for any any $R>0$, $f$ is uniformly continuous
on $[0,a]\times B_R\times B_R\times B_R$ and there exist positive constants $l_i$ ($i=1, 2, 3$) such that\vskip1mm
\begin{equation}\label{05}
\mu(f(t,U_1,U_2,U_3))\leq l_1\mu(U_1)+l_2\mu(U_2)+l_3\mu(U_3)
\end{equation}
\vskip2mm
\noindent for $\forall\;t\in [0,a]$ and bounded sets $U_1$, $U_2$, $U_3\subset E$, and
\vskip1mm\begin{equation}\label{06}
2a(l_1+aK_0l_2+aH_0l_3)<1.
\end{equation}\vskip2mm
\noindent Later, there are a large amount of authors studied ordinary
differential equations in Banach spaces similar to \eqref{04} by using
Sadovskii's fixed point theorem under the assumptions analogous to \eqref{05},
they also required that the constants satisfy a strong inequality
similar to \eqref{06}. For more details on this fact, please see Liu, Wu and Guo \cite{lwg04}, Liu et al. \cite{lg05} and the references therein.\vskip3mm

One can easily to discover that the inequality \eqref{03} and \eqref{06} are very strong restrictive
conditions, and they are difficult to be satisfied in applications. In
order to remove the strong restriction  on the constants in the conditions of noncompactness measure like \eqref{02} or \eqref{05}, Sun and
Zhang \cite{sz05} generalized the definition of condensing operator to
convex-power condensing operator. And based on the definition of
this new kind of operator, they established a new fixed point
theorem with respect to convex-power condensing operator which
generalizes the famous Sadovskii's fixed point theorem. As an application, the authors
investigated the existence of global mild solutions for the initial value problem (IVP) of evolution equations in the real Banach space $E$
\begin{equation}\label{07}
\left\{\begin{array}{ll}
u'(t)+Au(t)= f(t,u(t)),\quad t\in [0,a],\\[12pt]
u(0)=u_0,
  \end{array} \right.
\end{equation}
the authors assume the
nonlinear term $f$ is uniformly continuous on $[0,a]\times B_R$ and
satisfies a suitable noncompactness measure condition similar to
\eqref{02}. We should point out that the restriction condition similar to \eqref{03} has been deleted in \cite{sz05}.  Recently, Shi, Li and Sun \cite{sls11} developed the IVP \eqref{07} to the
case that the nonlinear term is $f(t,u(t),(Tu)(t), (Su)(t))$, and
obtained the existence of global mild solutions by using the new fixed point theorem with respect to
convex-power condensing operator established by Sun and
Zhang \cite{sz05}, but they also
require that the nonlinear term $f$ is uniformly continuous on
$[0,a]\times B_R\times B_R\times B_R$.\vskip3mm

During the past two decades, fractional differential equations and integral-differential equations have been proved to be valuable
tools in the investigation of many phenomena in engineering, physics, economy, chemistry,
aerodynamics and electrodynamics of complex medium. It has been found that the differential or integral-differential equations involving fractional derivatives in time are more realistic to describe many phenomena in practical cases than those of integer order in time. Particularly,
as the fractional order semilinear evolution equations are abstract formulations for many problems arising in engineering, physics, chemistry and economy, the theory of
fractional evolution equations has attracted increasing attention in recent years and it has developed into an important  branch of fractional calculus and fractional differential equations. For more details
about fractional calculus and fractional evolution equations, we refer to \cite{b01}, \cite{cp15}, \cite{cl14}, \cite{cl17}, \cite{ei04}, \cite{ei10}, \cite{gm98}, \cite{gl17}, \cite{lcl10},
 \cite{lpj12}, \cite{mpz15}, \cite{ss16}, \cite{wxl11}, \cite{wcx12}, \cite{wfz11}, \cite{wz11}, \cite{wzf13}, \cite{zj10}, \cite{zlw161}, \cite{zlw16} and the references therein. \vskip3mm

However, among the previous researches, most of researchers focus on the case that the differential operators in the main parts are
independent of time $t$, which means that the problems under considerations are autonomous.
On the other hand, we notice that when treating some parabolic evolution equations, it is usually assumed that the partial differential operators depend on time $t$ on account of this
class of operators appears frequently in the applications, for the details please see \cite{t97}.
As a result, it is significant and interesting to investigate fractional non-autonomous
evolution equations, i.e., the differential operators in the main parts of the considered
problems are dependent of time $t$. In fact, El-Borai \cite{ei04} investigated the existence and continuous dependence of fundamental solutions for a class of linear fractional non-autonomous evolution equations in 2004. In 2010, El-Borai, EI-Nadi and EI-Akabawy \cite{ei10} give some conditions to ensure the existence of resolvent operator for a class of fractional non-autonomous evolution equations with classical Cauchy initial condition.\vskip3mm

Motivated by the above mentioned aspects, in this paper we investigate the existence of mild solutions for the following Cauchy problem to nonlinear time fractional non-autonomous integro-differential evolution equation of mixed type via measure of noncompactness in Banach space $E$
\begin{equation}\label{11}
\left\{\begin{array}{ll}
  ^CD^{\alpha}_tu(t)+A(t)u(t)= f(t,u(t),(Tu)(t), (Su)(t)),\quad  t\in I, \\[12pt]
   u(0)=A^{-1}(0)u_0,
  \end{array} \right.
\end{equation}
where $ ^CD^{\alpha}_t$ is the standard Caputo's fractional time derivative of order $0<\alpha\leq 1$, $I=[0,a]$, $a>0$ is a constant, $A(t)$ is a family of closed linear operators defined on a dense domain $D(A)$ in Banach space $E$ into $E$ such
that $D(A)$ is independent of $t$,  $f:I\times E\times E\times E\rightarrow E$ is a Carath\'{e}odory type function, $u_0\in E$,
$T$ is the Volterra integral operator defined by \eqref{12} and
$S$ is the Fredholm integral operator defined by \eqref{13}.\vskip3mm

The motives and highlights in this article are as follows:
\begin{itemize}
\item [1.] We give the proper definition of mild solution for the Cauchy problem to  nonlinear time fractional non-autonomous integro-differential evolution equation of mixed type \eqref{11} by introducing three operators $\psi(t,s)$, $\varphi(t,\eta)$ and $U(t)$.\vskip4mm
\item [2.]  We observed that in \cite{g89}, \cite{ll81}, \cite{lg05}, \cite{lwg04}, \cite{sz05} and \cite{sls11},  the authors all demand that the
nonlinear term $f$ is uniformly continuous, this is a very strong
assumption. As a matter of fact, if $f(t,u)$ is Lipschitz continuous
on $I\times B_R$ with respect to $u$, then the condition \eqref{02} is
satisfied, but $f$ may not necessarily uniformly continuous on
$I\times B_R$. Therefore, in this work we deleted the assumption that $f$ is uniformly
continuous by using a new estimation technique of the measure of
noncompactness (see Lemma \ref{Le3}) established by Chen and Li in 2013 \cite{cl13}.\vskip4mm
\item [3.]  Just as we pointed out previously,  the inequality \eqref{06} is a very strong restrictive condition. How to get rid of the restriction on the constants in the conditions of noncompactness measure is a meaningful work. In this paper, we successfully used the new kind of fixed point theorem with respect to convex-power
condensing operator (see Lemma \ref{Le6}) established in 2005 to study time fractional non-autonomous evolution equations, and completely deleted the strong restriction on the constants in the conditions of noncompactness measure.
\end{itemize}\vskip3mm

The rest of this paper is organized as follows: in the sequel of Section 1, we give some general assumptions on the linear operator $-A(t)$, and also present the main results of this paper and its hypotheses. We provide in Section 2 some
definitions, notations and necessary preliminaries on fractional derivatives, Kuratowski measure of noncompactness and fixed point theorem with respect to convex-power condensing operator. In particular, the definition of mild solution for the Cauchy problem to time fractional non-autonomous integro-differential evolution equation of mixed type \eqref{11} and the properties about the operators $\psi(t,s)$, $\varphi(t,\eta)$ and $U(t)$ which defined in Definition \ref{De3} are also given. The proofs of the main theorems \ref{th1} and \ref{th2}, are given
in Section 3. In the final Section
4 we present an example of problem where the results of the previous
sections apply.

\vskip3mm
Let $\mathcal {L}(E)$ be the Banach space of all
linear and bounded operators in $E$ endowed with the topology defined by the operator norm.
Throughout this paper, we assume that the linear operator $-A(t)$ satisfies the following conditions:\vskip3mm
\begin{itemize}
\item [(A1)] For any $\lambda$ with $\textrm{Re}\lambda\geq 0$, the operator $\lambda I+A(t)$ exists a bounded inverse operator $\big[\lambda I+A(t)\big]^{-1}$ in $\mathcal {L}(E)$ and
    $$
    \Big\|\big[\lambda I+A(t)\big]^{-1}\Big\|\leq \frac{C}{|\lambda|+1},
    $$
where $C$ is a positive constant independent of both $t$ and $\lambda$;
\item [(A2)]  For any $t,\tau,s \in I$, there exists a constant $\gamma\in (0,1]$ such that
$$
\Big\|\big[A(t)-A(\tau)\big]A^{-1}(s)\Big\|\leq C|t-\tau|^\gamma,
$$
where the constants $\gamma$ and $C>0$ are independent of both $t$, $\tau$ and $s$.
\end{itemize}\vskip3mm
\noindent\textbf{Remark 1.}{\sl\quad From Henry \cite{he81}, Pazy \cite{pa83} and  Temam \cite{te97}, we know that the assumption (A1) means that for each $s\in J$, the operator $-A(s)$ generates an analytic
semigroup $e^{-tA(s)}$ ($t>0$), and there exists a positive constant $C$ independent of both $t$ and $s$ such that
$$
\Big\|A^n(s)e^{-tA(s)}\Big\|\leq \frac{C}{t^n},
$$
where $n=0$, $1$, $t>0$, $s\in J$.}
\vskip3mm
\noindent\textbf{Remark 2.}{\sl\quad In assumption (A1), if we choose $\lambda=0$ and $t=0$, then  there exists a positive constant $C$ independent of both $t$ and $\lambda$ such that
$$
\|A^{-1}(0)\|\leq C.
$$}
\noindent\textbf{Definition 3.}{\sl\quad  A function $\varphi:
[0,a]\times E\rightarrow E$ ia said to be Carath\'{e}odory continuous if
\begin{itemize}
\item [(i)] for all $ u\in E$,
 $\varphi(\cdot,u)$ is strongly measurable,\vskip 1mm
\item [(ii)] for a.e. $t\in [0,a]$, $\varphi(t,\cdot)$ is
 continuous.
 \end{itemize}
 }\vskip 3mm

In order to show the existence of mild solutions to the initial value problem for nonlinear time fractional non-autonomous integro-differential evolution equation of mixed type \eqref{11}, it is sufficient to impose some natural
growth conditions and noncompactness measure condition on the nonlinear function $f$.\vskip 3mm
\begin{itemize}
\item [(F1)]  For some $r>0$, there exist constants $0\leq\beta<\min\{\alpha,\gamma\}$,
$\rho_1>0$ and functions $\psi_r\in
L^{\frac{1}{\beta}}(J,\mathbb{R}^+)$ such that for a.e. $t\in I$ and
all $u\in E$ satisfying $\| u\|\leq r$, $$\|
f(t,u,Tu,Su)\|\leq\psi_r(t)\quad\textrm{and}\quad
\liminf\limits_{r\rightarrow +\infty}\frac{\|
\psi_r\|_{L^{\frac{1}{\beta}}[0,a]}}{r}:=\rho<+\infty.$$
\item [(F2)]  There exist positive constants $L_1$, $L_2$ and $L_3$ such that
for any bounded and countable sets $D_1$, $D_2$, $D_3\subset E$ and
a.e. $t\in I$,
$$\mu(f(t,D_1,D_2,D_3))\leq
L_1\mu(D_1)+L_2\mu(D_2)+L_3\mu(D_3).$$
\end{itemize}\vskip 3mm

For the sake of convenience, we denote by
$$
\delta_1=\Big(\frac{1-\beta}
{\alpha-\beta}\Big)^{1-\beta}+C\mathbf{B}(\alpha,\gamma)
a^{\gamma}\Big(\frac{1-\beta}{\alpha+\gamma-\beta}\Big)^{1-\beta}
$$
and
$$
\delta_2=1+Ca^\alpha\Big(\frac{1}{\alpha}+a^\gamma \mathbf{B}(\alpha,\gamma+1)\Big),
$$
where$$\mathbf{B}(\alpha,\gamma)=\int_0^1t^{\alpha-1}(1-t)^{\gamma-1}dt$$ is the Beta
function.\vskip 3mm

Our main results are as follows:\vskip 2mm
\begin{theorem}\label{th1}
Assume that the nonlinear function $f:I\times E\times E\times E\rightarrow E$ is Carath\'{e}odory continuous. If the assumptions (F1) and (F2) are satisfied,
then problem \eqref{11} exists at least one mild solution in
$C(I,E)$ provided that
\begin{equation} \label{1}
C\rho a^{\alpha-\beta}\delta_1< 1.
\end{equation}
\end{theorem}\vskip 3mm

\begin{theorem}\label{th2}
Assume that the nonlinear function $f:I\times E\times E\times E\rightarrow E$ is Carath\'{e}odory continuous. If the assumption (F2) and the following assumptions
\begin{itemize}
\item [(F1)$^{\ast}$]  There exist a functions $\phi\in
L^{\frac{1}{\beta}}(I,\mathbb{R}^+)$  for $0\leq\beta<\min\{\alpha,\gamma\}$ and a nondecreasing continuous function $\Phi:\mathbb{R}^+\rightarrow
\mathbb{R}^+$ such that$$\|
f(t,u,Tu,Su)\|\leq\phi(t)\Phi(\|u\|)$$ for a.e. $t\in I$ and
all $u\in E$,
\end{itemize}
hold, then problem \eqref{11} exists at least one mild solution in
$C(I,E)$ provided that  there exists a positive constant $R$ such that
\begin{equation} \label{2}
\delta_1\Phi(R)a^{\alpha-\beta}\|
\phi \|_{L^{\frac{1}{\beta}}[0,a]}\leq \frac{R}{C}-\delta_2\|u_0\|.
\end{equation}
\end{theorem}

\vskip 3mm
From Theorem \ref{th2}, we can obtain the following result.

\begin{corollary}\label{Co1}
Assume that the nonlinear function $f:I\times E\times E\times E\rightarrow E$ is Carath\'{e}odory continuous. If the assumptions (F1)$^{\ast}$ and (F2) are satisfied, then problem \eqref{11} exists at least one mild solution in
$C(I,E)$ provided that
\begin{equation} \label{3}
\liminf\limits_{r\rightarrow+\infty}\frac{\Phi(r)}{r}<\frac{1}{C\delta_1a^{\alpha-\beta}\|
\phi \|_{L^{\frac{1}{\beta}}[0,a]}}.
\end{equation}
\end{corollary}
\vskip3mm
\noindent\textbf{Remark 4.}{\sl\quad In Theorems \ref{th1}, \ref{th2} and Corollary \ref{Co1}, we did not make any restrictions on the constants $L_1$, $L_2$ and $L_3$ in the assumption (F2). Noticing that in \cite{ll81}, the constant $l$ should satisfy the restriction condition \eqref{03} and in \cite{g89}, the constants $l_1$, $l_2$, $l_3$ should satisfy the restriction condition \eqref{06}. In this paper, we deleted the restrictions on the constants $L_1$, $L_2$ and $L_3$, this is a huge improvement. Furthermore, in this paper we only assume that the nonlinear function $f$ is Carath\'{e}odory continuous, which is weaker than $f$ is uniformly continuous required in \cite{g89}, \cite{ll81}, \cite{lg05}, \cite{lwg04}, \cite{sz05} and \cite{sls11}.}

\vskip 15mm
\section{Preliminaries}
In this section, we introduce some notations, definitions, and preliminary facts  including fractional derivatives and integrals, Kuratowski measure of noncompactness, fixed point theorem with respect to convex-power condensing operator and the operators $\psi(t,s)$, $\varphi(t,\eta)$ and $U(t)$, which are used throughout
this paper.\vskip 3mm

Throughout this paper, we set $I=[0,a]$ denotes a compact interval in $\mathbb{R}$, where $a>0$ is a constant. Let $E$ be a Banach space with norm $\|\cdot\|$. We denote by $C(I,E)$ the Banach space of
all continuous functions from interval $I$ into $E$ equipped with the supremum norm
$$\| u\|_{C}=\sup\{\|u(t)\|,\quad t\in I\} \qquad \forall\; u\in C(I,E),$$
and by $\mathcal {L}(E)$ the Banach space of all
linear and bounded operators in $E$ endowed with the topology defined by the operator norm. Let $L^1(I,E)$ be the Banach space of all
$E$-value Bochner integrable functions defined on $I$ with the norm
$\|
u\|_1=\int_0^{a}\| u(t)\| dt$. For any $r>0$, denote $
\Omega_r=\{u\in C(I,E): \|u(t)\|\leq r,\;t\in I\}$,
then $\Omega_r$ is a closed ball in $C(I,E)$ with center $\theta$ and radius $r$.\vskip 3mm

At first, we recall the definition of the Riemann-Liouville integral and Caputo derivative of fractional order.\vskip 3mm

\begin{definition}\label{De1}(\cite{co71}). The fractional integral of
order $\alpha>0$ with the lower limit zero for a function $f\in L^1([0,+\infty),\mathbb{R})$
is defined as
$$I^\alpha_t f(t)=\frac{1}{\Gamma(\alpha)}\int_0^t (t-s)
^{\alpha-1}f(s)ds.$$ Here and elsewhere
$$
\Gamma(\alpha)=\int_0^{+\infty}t^{\alpha-1}e^{-t}dt
$$ denotes the Gamma
function.
\end{definition}\vskip 3mm

\begin{definition}\label{De2}(\cite{co71}). The Caputo fractional
derivative of order $\alpha$ with the lower limit zero for a function
$f:[0,+\infty)\rightarrow\mathbb{R}$, which is at least $n$-times differentiable can be defined as
$$^CD^{\,\alpha}_t f(t)=\frac{1}{\Gamma(n-\alpha)}\int_0^t (t-s)
^{n-\alpha-1}f^{(n)}(s)ds=I_t^{n-\alpha}f^{(n)}(t),$$ where $n-1<\alpha<n$, $n\in \mathbb{N}$.
\end{definition}
\vskip3mm
\begin{lemma}\label{Le0} \textrm{(Bochner's Theorem)}
A measurable function $f$ maps from $[0,+\infty)$ to $E$ is Bochner integrable if
$\|f\|$ is Lebesgue integrable.
\end{lemma}\vskip3mm
\begin{remark}\label{Re1}
(i) $I^{\alpha_1}_tI^{\alpha_2}_t=I^{\alpha_1+\alpha_2}_t$ for any $\alpha_1$, $\alpha_2>0$;\\[6pt]
(ii) If $f$ is an abstract function with values in $E$, then the integrals which appear in
Definitions \ref{De1} and \ref{De2} are taken in Bochner's sense.
\end{remark}
\vskip3mm

By the above discussion and \cite[Theorem 2.6]{ei04},
we can get the definition of mild solutions for problem \eqref{11}.\vskip 3mm

\begin{definition}\label{De3} A function $u\in C(I,E)$ is said
to be a mild solution of the initial value problem \eqref{11} if it satisfies
\begin{equation} \label{24}\begin{split}
u(t)
&=A^{-1}(0)u_0+\int_0^t\psi(t-\eta,\eta)U(\eta)u_0d\eta\\
&\quad+\int_0^t\psi(t-\eta,\eta)
f(\eta,u(\eta),(Tu)(\eta),(Su)(\eta))d\eta
\\
&\quad+\int_0^t\int_0^\eta\psi(t-\eta,\eta)\varphi(\eta,s)
f(s,u(s),(Tu)(s),(Su)(s))dsd\eta,\end{split}
\end{equation}
where the operators  $\psi(t,s)$, $\varphi(t,\eta)$ and $U(t)$ are defined by
\begin{equation} \label{21}
\psi(t,s)=\alpha\int_0^\infty\theta t^{\alpha-1}\xi_\alpha(\theta)e^{-t^\alpha \theta A(s)}d\theta,
\end{equation}
 \begin{equation} \label{22}
\varphi(t,\eta)=\sum
_{k=1}^{\infty}\varphi_k(t,\eta)
\end{equation}and
\begin{equation} \label{23}
U(t)=-A(t)A^{-1}(0)-\int_0^t\varphi(t,s)A(s)A^{-1}(0)ds,
\end{equation}
$\xi_\alpha$ is a probability density function defined on $[0,+\infty)$ such that it's Laplace transform is
given by
$$
\int_0^\infty e^{-\theta x}\xi_\alpha(\theta)d\theta=\sum
_{i=0}^{\infty}\frac{(-x)^i}{\Gamma(1+\alpha i)},\quad 0<\alpha\leq 1,\; x>0,
$$\vskip1mm
$$
\varphi_1(t,\eta)=\big[A(t)-A(\eta)\big]\psi(t-\eta,\eta),
$$
and
$$
\varphi_{k+1}(t,\eta)=\int_\eta^t\varphi_{k}(t,s)\varphi_1(s,\eta)d\eta,\quad k=1,2,\cdots.
$$
\end{definition}\vskip3mm

The following properties about the operators  $\psi(t,s)$, $\varphi(t,\eta)$ and $U(t)$ will be needed in our argument.\vskip 3mm

\begin{lemma}\label{Le1}(\cite{ei04}) The operator-valued functions $\psi(t-\eta,\eta)$ and $A(t)\psi(t-\eta,\eta)$ are continuous
in uniform topology about the variables $t$ and $\eta$, where $t\in I$, $0\leq\eta\leq t-\epsilon$  for any $\epsilon>0$, and
\begin{equation} \label{25}
\|\psi(t-\eta,\eta)\|\leq C(t-\eta)^{\alpha-1},
\end{equation}
where $C$ is a positive constant independent of both $t$ and $\eta$. Furthermore,
\begin{equation} \label{26}
\|\varphi(t,\eta)\|\leq C(t-\eta)^{\gamma-1}
\end{equation}
and
\begin{equation} \label{27}
\|U(t)\|\leq C(1+t^\gamma).
\end{equation}
\end{lemma}\vskip3mm

By Lemma \ref{Le1} and appropriate calculation, we can obtain the following result.
\vskip3mm
\begin{lemma}\label{Le21} The integral $\int_0^t\psi(t-\eta,\eta)U(\eta)d\eta$ is uniformly continuous in the operator norm $\mathcal {L}(E)$ for any $t\in I$, and
\begin{equation} \label{28}
\Big\|\int_0^t\psi(t-\eta,\eta)U(\eta)d\eta\Big\|\leq C^2t^\alpha\Big(\frac{1}{\alpha}+t^\gamma \mathbf{B}(\alpha,\gamma+1)\Big),\qquad \forall\;t\in I.
\end{equation}
\end{lemma}\vskip3mm

Using the properties of Riemann-Liouville integral of fractional order and proper integral transformation, we can obtain the following result, which will be used in the proof of our main results.\vskip3mm

\begin{lemma}\label{Le22}For any $t\in I$ and $g\in L^1[0,a]$, we have
\begin{equation} \label{29}
\int_0^t\int_0^\eta(t-\eta)^{\alpha-1}(\eta-s)^{\gamma-1}g(s)dsd\eta=\mathbf{B}(\alpha,\gamma) \int_0^t(t-\eta)^{\alpha+\gamma-1}g(\eta)d\eta.
\end{equation}
\end{lemma}\vskip3mm

Next, we introduce the definition for Kuratowski measure of noncompactness, which will be used in the proof of our main results.\vskip3mm
\begin{definition}\label{De4} (\cite{bg80}) The Kuratowski measure of noncompactness $\mu(\cdot)$ defined on bounded set $S$ of Banach space $E$ is
$$
\mu(S):=\inf\{\delta>0\,:\,S=\cup_{k=1}^{n}S_k\; \textrm{and}\; \textrm{diam}(S_k)\leq \delta\,\textrm{ for }\,k=1,2,\cdots,n\}.
$$
\end{definition}\vskip3mm

The following properties about the Kuratowski measure of noncompactness are well known.
\vskip3mm
\begin{lemma}\label{Le2}(\cite{d85})  Let $E$ be a Banach space and $U$, $V\subset E$ be bounded. The following properties are satisfied :
\begin{itemize}
\item [\textrm{(1)}] $\mu(U)\leq\mu(V)$ if $U\subset V$;\vskip2mm
\item [\textrm{(2)}] $\mu(U)=\mu(\overline{U})=\mu(\textrm{conv}\; U)$, where \textrm{conv} $U$ means the convex hull of $U$;\vskip2mm
\item [\textrm{(3)}] $\mu(U)=0$ if and only if $\overline{U}$ is compact, where $\overline{U}$ means the closure hull of $U$;\vskip2mm
\item [\textrm{(4)}] $\mu(\lambda U)=|\lambda|\mu(U)$, where $\lambda\in \mathbb{R}$;\vskip2mm
\item [\textrm{(5)}] $\mu(U\cup V)=\max\{\mu(U),\mu(V)\}$;\vskip2mm
\item [\textrm{(6)}] $\mu(U+V)\leq\mu(U)+\mu(V)$, where $U+V=\{x\mid x=y+z, y\in U, z\in V\}$;\vskip2mm
\item [\textrm{(7)}] $\mu(U+x)=\mu(U)$, for any $x\in E$;\vskip2mm
\item [\textrm{(8)}] If the map $Q:D(Q)\subset E\rightarrow X$ is Lipschitz continuous with constant $k$, then $\mu(Q(S))\leq k\mu(S)$ for any bounded subset $S\subset \mathcal {D}(Q)$, where $X$ is another Banach space.
\end{itemize}\end{lemma}
\vskip3mm

In this article, we denote by $\mu(\cdot)$ and $\mu_C(\cdot)$ the Kuratowski measure of noncompactness on the bounded set of $E$ and $C(I,E)$ respectively. For any $D\subset C(I,E)$ and $t\in I$, set $D(t)=\{u(t)\mid u\in
D\}$ then $D(t)\subset E$. If $D\subset C(I,E)$ is bounded,
then $D(t)$ is bounded in $E$ and $\mu(D(t))\leq \mu_C(D)$. For more details about the properties of the Kuratowski measure of noncompactness, we refer to the monographs \cite{bg80} and \cite{d85}.\vskip3mm

The following lemmas are needed in our argument.\vskip3mm
\begin{lemma}\label{Le5}(\cite{bg80}).
 Let $E$ be a Banach space, and let $D\subset C(I,E)$ be bounded and equicontinuous.
 Then $\mu(D(t))$ is continuous on $[0,a]$, and
$ \mu_C(D)=\max\limits_{t\in [0,a]}\mu(D(t))$.
\end{lemma}\vskip3mm
\begin{lemma}\label{Le3} (\cite{cl13}).
Let $E$ be a Banach space,
and let $D\subset E$ be bounded. Then there exists a countable set
$D_0\subset D$, such that $\mu(D)\leq 2\mu(D_0)$.
\end{lemma}\vskip3mm

\begin{lemma}\label{Le4}(\cite{h83}).
Let $E$ be a Banach space. If $D=\{u_n\}_{n=1}^\infty\subset C([0,a],E)$ is a countable set and there exists a function $m\in L^1([0,a],\mathbb{R}^+)$ such that for every $n\in \mathbb{N}$
$$
\|u_n(t)\|\leq m(t),\qquad \textrm{a.e.}\; t\in [0,a].
$$
Then
$\mu(D(t))$ is Lebesgue integral on $[0,a]$, and
$$
\mu\Big(\Big\{\int_{0}^{a} u_n(t)dt\mid n\in \mathbb{N}\Big\}\Big)\leq
2\int_{0}^{a}\mu(D(t))dt.
$$
\end{lemma}
\vskip3mm

The following fixed point theorem with respect to convex-power condensing
operator which introduced by Sun and Zhang \cite{sz05} plays a key role in the proof of our main results.\vskip3mm

\begin{definition}\label{De5} Let $E$ be a real Banach space. If $\mathcal {A}:E\rightarrow E$ is a continuous and bounded operator, there exist $u_0\in E$ and a positive integer $n_0$ such that for
any bounded and nonprecompact subset $S\subset E$,\\
\begin{equation} \label{210}\mu(\mathcal {A}^{(n_0,u_0)}(S))< \mu(S),\end{equation}
where
$$\mathcal {A}^{(1,u_0)}(S)\equiv \mathcal {A}(S),\quad \mathcal {A}^{(n,u_0)}(S)= \mathcal {A}(\overline{co}\{Q^{(n-1,u_0)}(S),u_0\}),\quad n=2,3,\cdots.$$
Then we call $\mathcal {A}$ a convex-power condensing operator about $u_0$ and $n_0$.
\end{definition}\vskip3mm

\begin{lemma}\label{Le6}(Fixed point theorem with respect to convex-power condensing operator, \cite{sz05})  Let $E$ be a real Banach space, and let
$D\subset E$
be a bounded, closed and convex set in $E$. If there exist $u_0\in D$ and a positive integer $n_0$ such that $\mathcal {A}:D\rightarrow D$ be a convex-power condensing operator about $u_0$ and $n_0$, then the operator $\mathcal {A}$ exists at least one fixed point in $D$.
\end{lemma}
\vskip3mm

\begin{remark}\label{re1}
If $n_0=1$ in \eqref{210}, then fixed point theorem with respect to convex-power condensing operator (see Lemma \ref{Le6}) will degrade into the famous Sadoveskii's fixed point theorem (see Theorem 2). Noticed that Lemma \ref{Le6} requires the operator $\mathcal {A}$ is neither condensing nor completely continuous. Therefore, fixed point theorem with respect to convex-power condensing operator is the
generalization of the well-known Sadoveskii's fixed point theorem.
\end{remark}
\vskip15mm

\section{Proof of the main results}
In this section, we give the proofs of Theorem \ref{th1} and Theorem \ref{th2}.
\vskip2mm
\noindent\textbf{Proof of Theorem \ref{th1}.}\quad Define an operator $\mathcal {Q}$ on the space of continuous functions $C(I,E)$ as follows
\begin{equation} \label{31}\begin{split}
(\mathcal {Q}u)(t)
&=A^{-1}(0)u_0+\int_0^t\psi(t-\eta,\eta)U(\eta)u_0d\eta\\[6pt]
&\quad+\int_0^t\psi(t-\eta,\eta)
f(\eta,u(\eta),(Tu)(\eta),(Su)(\eta))d\eta
\\[6pt]
&\quad+\int_0^t\int_0^\eta\psi(t-\eta,\eta)\varphi(\eta,s)
f(s,u(s),(Tu)(s),(Su)(s))dsd\eta.\end{split}
\end{equation}
By direct calculation and the properties about the operators  $\psi(t,s)$, $\varphi(t,\eta)$ and $U(t)$, it is easy to know that the operator $\mathcal {Q}$ maps $C(I,E)$ to $C(I,E)$, and it is well defined. From
Definition \ref{De3}, one can easily to verify that the mild solution of initial value problem \eqref{11} is equivalent to the fixed point of the operator $\mathcal {Q}$ defined by \eqref{31}. In
what follows, we will prove that the operator $\mathcal {Q}$  has at least one fixed point by applying the fixed point theorem with respect to convex-power condensing operator.
\vskip3mm

Firstly, we prove that there exists a positive constant $R$ such
that the operator $\mathcal {Q}$ defined by \eqref{31} maps the set $\Omega_{R}$ to $\Omega_{R}$. If this is not true, then there would exist $t_r\in I$ and $u_r\in \Omega_{r}$
such that $\| (\mathcal {Q}u_r)(t_r)\|>r$ for
each $r>0$. Combining with Lemmas
\ref{Le1}-\ref{Le22}, the assumption (F1) and H\"{o}lder inequality, we get that
\begin{eqnarray*}
r&<&\| (\mathcal {Q}u_r)(t_r)\|\leq \|A^{-1}(0)u_0\|+\Big\|\int_0^{t_r}
\psi({t_r}-\eta,\eta)U(\eta)u_0d\eta\Big\|\\& &+\Big\|\int_0^{t_r}\psi({t_r}-\eta,\eta)
f(\eta,u_r(\eta),(Tu_r)(\eta),(Su_r)(\eta))d\eta\Big\|\\[6pt]& &
+\Big\|\int_0^{t_r}\int_0^\eta\psi({t_r}-\eta,\eta)\varphi(\eta,s)
f(s,u_r(s),(Tu_r)(s),(Su_r)(s))dsd\eta\Big\|\\ &\leq & C\|u_0\|+C^2\int_0^{t_r}({t_r}-\eta)^{\alpha-1}(1+\eta^\gamma)\|u_0\|d\eta
\qquad\qquad\qquad\qquad\qquad\qquad\qquad\qquad\qquad\qquad\qquad\qquad\qquad\;
\end{eqnarray*}
\begin{equation} \label{32}\begin{split} 
&\quad +C\int_0^{t_r}({t_r}-\eta)^{\alpha-1}\psi_r(\eta)d\eta
+C^2\int_0^{t_r}\int_0^\eta ({t_r}-\eta)^{\alpha-1}(\eta-s)^{\gamma-1}\psi_r(s)dsd\eta\\
&\leq C\|u_0\|
+C^2\|u_0\|({t_r})^\alpha\Big(\frac{1}{\alpha}+({t_r})^\gamma \mathbf{B}(\alpha,\gamma+1)\Big)\\
&\quad +C\int_0^{t_r}({t_r}-\eta)^{\alpha-1}\psi_r(\eta)d\eta
+C^2\mathbf{B}(\alpha,\gamma)\int_0^{t_r} ({t_r}-\eta)^{\alpha+\gamma-1}\psi_r(\eta)d\eta\\
&\leq C\|u_0\|
+C^2\|u_0\|a^\alpha\Big(\frac{1}{\alpha}+a^\gamma \mathbf{B}(\alpha,\gamma+1)\Big)\\
&\quad +C\Big(\int_0^{t_r}({t_r}-\eta)^{\frac{\alpha-1}{1-\beta}}d\eta\Big)^{1-\alpha_1} \Big(\int_0^{t_r}\psi_r^{\frac{1}{\beta}}(\eta)d\eta\Big)^{\beta}\\
&\quad+C^2\mathbf{B}(\alpha,\gamma)\Big(\int_0^{t_r}({t_r}-\eta)^{\frac{\alpha+\gamma-1}
{1-\beta}}d\eta\Big)^{1-\beta} \Big(\int_0^{t_r}\psi_r^{\frac{1}{\beta}}(\eta)d\eta\Big)^{\beta}\\
&\leq C\|u_0\|
+C^2\|u_0\|a^\alpha\Big(\frac{1}{\alpha}+a^\gamma \mathbf{B}(\alpha,\gamma+1)\Big)\\
&\quad+ Ca^{\alpha-\beta}\Big(\frac{1-\beta}{\alpha-\beta}\Big)^{1-\beta}\|
\psi_r\|_{L^{\frac{1}{\beta}}[0,a]}\\
&\quad+C^2\mathbf{B}(\alpha,\gamma)a^{\alpha+\gamma-\beta}
\Big(\frac{1-\beta}{\alpha+\gamma-\beta}\Big)^{1-\beta}\|
\psi_r\|_{L^{\frac{1}{\beta}}[0,a]}.\end{split}
\end{equation}
Dividing both side of \eqref{32} by $r$ and taking the lower limit as
$r\rightarrow +\infty$, combined with the assumption \eqref{1} we get that
\begin{equation} \label{33}\begin{split}
1&\leq C\rho a^{\alpha-\beta}\Big[\Big(\frac{1-\beta}
{\alpha-\beta}\Big)^{1-\beta}+C\mathbf{B}(\alpha,\gamma)
a^{\gamma}\Big(\frac{1-\beta}{\alpha+\gamma-\beta}\Big)^{1-\beta}\Big]\\[8pt]&=C\rho a^{\alpha-\beta}\delta_1<1.
\end{split}
\end{equation}
One can easily to see that \eqref{33} is a contradiction. Therefore, we have proved that $\mathcal {Q}:\Omega_{R}\to \Omega_{R}$.
\vskip3mm

Secondly, we prove that the operator $\mathcal {Q}:\Omega_{R}\to \Omega_{R}$ is continuous. To this
end, let $\{u_n\}_{n=1}^\infty\subset \Omega_{R}$ be a sequence such
that $\lim\limits_{n\rightarrow +\infty}u_n=u$ in $\Omega_{R}$. By
the Carath\'{e}odory continuity of the nonlinear function $f$, we get that
\begin{equation} \label{34}\lim\limits_{n\rightarrow
+\infty}\|f(t,u_n(t),(Tu_n)(t),(Su_n)(t))-f(t,u(t),(Tu)(t),(Su)(t))\|
=0
\end{equation}for a.e.
$t\in I$.
By \eqref{31} and Lemmas
\ref{Le1}-\ref{Le22} combined with the similar calculus method with which used in \eqref{32},  we have

\begin{equation} \label{36}\begin{split}
 & \| (\mathcal {Q}u_n)(t)-(\mathcal {Q}u)(t)\|\\&\leq
\Big\|\int_0^t\psi(t-\eta,\eta)[
f(\eta,u_n(\eta),(Tu_n)(\eta),(Su_n)(\eta))\\[6pt]&\quad
-f(t,u(\eta),(Tu)(\eta),(Su)(\eta))]d\eta\Big\|\\
&\quad +\Big\|\int_0^t\int_0^\eta\psi(t-\eta,\eta)\varphi(\eta,s)[
f(s,u_n(s),(Tu_n)(s),(Su_n)(s))\\[6pt]&\quad
-f(t,u(s),(Tu)(s),(Su)(s))]dsd\eta\Big\|\\
&\leq
C\int_0^{t}({t}-\eta)^{\alpha-1}\|f(\eta,u_n(\eta),(Tu_n)(\eta),(Su_n)(\eta))\\[6pt]&\quad
-f(t,u(\eta),(Tu)(\eta),(Su)(\eta))\|d\eta\\
&\quad + C^2\int_0^{t}\int_0^\eta ({t}-\eta)^{\alpha-1}(\eta-s)^{\gamma-1}\|f(s,u_n(s),(Tu_n)(s),(Su_n)(s))\\[6pt]&\quad
-f(t,u(s),(Tu)(s),(Su)(s))\|ds d\eta,\qquad \forall\; t\in I.
\end{split}
\end{equation}
By the assumption (F1), we know that for every $t\in I$,

\begin{equation} \label{37}\begin{split}
({t}-\eta)^{\alpha-1}&\|f(\eta,u_n(\eta),(Tu_n)(\eta),(Su_n)(\eta))\\[6pt]&\quad
-f(t,u(\eta),(Tu)(\eta),(Su)(\eta))\|\\[6pt]
&\leq 2({t}-\eta)^{\alpha-1}\psi_R(\eta)
\end{split}
\end{equation}
for a.e. $\eta \in[0,t]$.
By again the assumption (F1) combined with Lemma \ref{Le22}, we get that for every $t\in I$, $0\leq\eta\leq t$ and a.e. $s\in[0,\eta]$,
\begin{equation} \label{38}\begin{split}
& \int_0^{t}\int_0^\eta ({t}-\eta)^{\alpha-1}(\eta-s)^{\gamma-1}\|f(s,u_n(s),(Tu_n)(s),(Su_n)(s))\\[6pt]&\quad
-f(s,u(s),(Tu)(s),(Su)(s))\|ds d\eta\\[6pt]
&\leq 2\int_0^{t}\int_0^\eta ({t}-\eta)^{\alpha-1}(\eta-s)^{\gamma-1} \psi_R(s)ds d\eta\\[6pt]
&= 2\mathbf{B}(\alpha,\gamma)\int_0^{t}({t}-\eta)^{\alpha+\gamma-1}\psi_R(\eta)d\eta.
\end{split}
\end{equation}
From the fact that the functions $\eta\rightarrow
2({t}-\eta)^{\alpha-1}\psi_R(\eta)$ and $\eta\rightarrow
2\mathbf{B}(\alpha,\gamma)({t}-\eta)^{\alpha+\gamma-1}\psi_R(\eta)$ are Lebesgue integrable for a.e. $\eta\in[0,t]$ and every
$t\in I$, combined with \eqref{34}, \eqref{36}, \eqref{37}, \eqref{38} and the Lebesgue dominated
convergence theorem, we know that
$$
\| (\mathcal {Q}u_n)(t)-(\mathcal {Q}u)(t)\|\to 0\quad \textrm{as}\quad n\rightarrow \infty
$$
 for any $t\in I$. Therefore, we get that
$$\|\mathcal {Q}u_n-\mathcal {Q}u\|_C\rightarrow 0\quad (n\rightarrow \infty),$$
which means that $\mathcal {Q}:\Omega_{R}\to \Omega_{R}$ is a continuous operator.\vskip3mm

Now, we are in the position to demonstrate that $\mathcal {Q}:\Omega_{R}
\rightarrow\Omega_{R}$ is an equicontinuous operator. For any $u\in\Omega_{R}$
and $0\leq t'< t''\leq a$, by \eqref{31} and the assumption (F1), we know that

\begin{eqnarray*}
& &\|(\mathcal {Q}u)(t'')-(\mathcal {Q}u)(t')\|\\&\leq& \Big\|\int_{t'}^{t''}\psi(t''-\eta,\eta)U(\eta)u_0d\eta\Big\|\\
& &+\Big\|\int_0^{t'}[\psi(t''-\eta,\eta)-\psi(t'-\eta,\eta)]U(\eta)u_0d\eta\Big\|\\
& &+\Big\|\int_{t'}^{t''}\psi(t''-\eta,\eta)
f(\eta,u(\eta),(Tu)(\eta),(Su)(\eta))d\eta\Big\|\\
& &+\Big\|\int_{0}^{t'}[\psi(t''-\eta,\eta)-\psi(t'-\eta,\eta)]
f(\eta,u(\eta),(Tu)(\eta),(Su)(\eta))d\eta\Big\|\\
& &+\Big\|\int_{t'}^{t''}\int_0^\eta\psi(t''-\eta,\eta)\varphi(\eta,s)
f(s,u(s),(Tu)(s),(Su)(s))dsd\eta\Big\|\\
& &+\Big\|\int_{0}^{t'}\int_0^\eta[\psi(t''-\eta,\eta)-\psi(t'-\eta,\eta)]\varphi(\eta,s)\\
& &\quad\cdot
f(s,u(s),(Tu)(s),(Su)(s))dsd\eta\Big\|\\[6pt]
&\leq&J_1+J_2+J_3+J_4+J_5+J_6,
\end{eqnarray*}
where
\begin{eqnarray*}
\begin{array}{lc}
J_1=\displaystyle\int_{t'}^{t''}\|\psi(t''-\eta,\eta)U(\eta)\|\cdot\|u_0\|d\eta,\\[12pt]
J_2=\displaystyle\int_0^{t'}\|[\psi(t''-\eta,\eta)-\psi(t'-\eta,\eta)]U(\eta)\|\cdot\|u_0\|d\eta,
\end{array}
\end{eqnarray*}
\begin{eqnarray*}
\begin{array}{lc}
J_3=\displaystyle\int_{t'}^{t''}\|\psi(t''-\eta,\eta)\|\psi_R(\eta)d\eta,\\[12pt]
J_4=\displaystyle\int_{0}^{t'}\|\psi(t''-\eta,\eta)-\psi(t'-\eta,\eta)\|\psi_R(\eta)d\eta,
\\[12pt]
J_5=\displaystyle\int_{t'}^{t''}\int_0^\eta\|\psi(t''-\eta,\eta)\varphi(\eta,s)\|\psi_R(s)dsd\eta,
\\[12pt]
J_6=\displaystyle\int_{0}^{t'}\int_0^\eta\|[\psi(t''-\eta,\eta)-
\psi(t'-\eta,\eta)]\varphi(\eta,s)\|\psi_R(s)dsd\eta.
\end{array}
\end{eqnarray*}
Therefore, we only need to prove that $J_k\to 0$ independently of $u\in \Omega_{R}$ as $t''-t'\to 0$ for $k=1,\,2,\,3,\,4,\,5,\,6$.\vskip3mm
\noindent For $J_1$, by Lemma \ref{Le1} we know that
$$
J_1\leq C^2\|u_0\|\int_{t'}^{t''}(t''-\eta)^{\alpha-1}(1+\eta^\gamma)d\eta\rightarrow0 \quad \textrm{as}\quad t''-t'\rightarrow 0.
$$\vskip3mm
\noindent
For $t'=0$ and $0<t''\leq a$, it is easy to see that $J_2=0$. For $t'>0$ and $\epsilon>0$ small enough, by Lemma \ref{Le1} and the fact that operator-valued function $\psi(t-\eta,\eta)$ is continuous
in uniform topology about the variables $t$ and $\eta$ for $0\leq t\leq a$ and $0\leq\eta\leq t-\epsilon$, we have
\begin{eqnarray*}
J_2&\leq&
\sup\limits_{\eta\in
[0,t'-\epsilon]}\|\psi(t''-\eta,\eta)-\psi(t'-\eta,\eta)\|\cdot C\|u_0\|
\int_{0}^{t'-\epsilon}
(1+\eta^\gamma)d\eta\\[6pt]
& &+C^2\|u_0\|\int_{t'-\epsilon}^{t'}
[(t''-\eta)^{\alpha-1}+(t'-\eta)^{\alpha-1}](1+\eta^\gamma)d\eta\\[6pt]
&\rightarrow&0\quad \textrm{as}\quad t''-t'\rightarrow 0\quad \textrm{and}\quad \epsilon\rightarrow 0.
\end{eqnarray*}\vskip3mm
\noindent
For $J_3$, by Lemma \ref{Le1}, the assumption (F1) and H\"{o}lder inequality, we get that
\begin{eqnarray*}
J_3&\leq&
C\int_{t'}^{t''}
(t''-\eta)^{\alpha-1}\psi_R(\eta)d\eta\\[6pt]
&\leq&C\Big(\int_{t'}^{t''}({t''}-\eta)^{\frac{\alpha-1}{1-\beta}}d\eta\Big)^{1-\beta} \Big(\int_{t'}^{t''}\psi_R^{\frac{1}{\beta}}(\eta)d\eta\Big)^{\beta}\\[6pt]
&\leq&C\Big(\frac{1-\beta}{\alpha-\beta}\Big)^{1-\beta}\|
\psi_R\|_{L^{\frac{1}{\beta}}[0,a]}(t''-t')^{\alpha-\beta}\\[6pt]
&\rightarrow&0\quad \textrm{as}\quad t''-t'\rightarrow 0.
\end{eqnarray*}\vskip3mm
\noindent
For $t'=0$ and $0<t''\leq a$, it is easy to see that $J_4=0$. For $t'>0$ and $\epsilon>0$ small enough, by Lemma \ref{Le1} and the fact that operator-valued function $\psi(t-\eta,\eta)$ is continuous
in uniform topology about the variables $t$ and $\eta$ for $0\leq t\leq a$ and $0\leq\eta\leq t-\epsilon$, we know that

\begin{eqnarray*}
J_4&\leq&
(t'-\epsilon)^{1-\beta}\|
\psi_R\|_{L^{\frac{1}{\beta}}[0,a]}\sup\limits_{\eta\in
[0,t'-\epsilon]}\|\psi(t''-\eta,\eta)-\psi(t'-\eta,\eta)\|\\[6pt]
& &+C\int_{t'-\epsilon}^{t'}
[(t''-\eta)^{\alpha-1}+(t'-\eta)^{\alpha-1}]\psi_R(\eta)d\eta
\\[6pt]
&\rightarrow&0\quad \textrm{as}\quad t''-t'\rightarrow 0\quad \textrm{and}\quad \epsilon\rightarrow 0.
\end{eqnarray*}\vskip5mm
\noindent
For $J_5$, by Lemma \ref{Le1}, the assumption (F1) and the fact that the function $\eta\to(t''-\eta)^{\alpha-1}I_\eta^\gamma \psi_R(\eta)$ is Lebesgue integrable, we have
\begin{eqnarray*}
J_5&\leq&
C^2\int_{t'}^{t''}\int_{0}^{\eta}
(t''-\eta)^{\alpha-1}(\eta-s)^{\gamma-1}\psi_R(s)dsd\eta
\\[6pt]&\leq&C^2\Gamma(\gamma)\int_{t'}^{t''}
(t''-\eta)^{\alpha-1}I_\eta^\gamma \psi_R(\eta)d\eta
\\[6pt]&\rightarrow&0\quad \textrm{as}\quad t''-t'\rightarrow 0.
\end{eqnarray*}\vskip3mm
\noindent
For $t'=0$ and $0<t''\leq a$, it is easy to see that $J_6=0$. For $t'>0$ and $\epsilon>0$ small enough, by Lemma \ref{Le1}, the assumption (F1), the facts that the functions $\eta\to(t''-\eta)^{\alpha-1}I_\eta^\gamma \psi_R(\eta)$ and $\eta\to(t'-\eta)^{\alpha-1}I_\eta^\gamma \psi_R(\eta)$ are Lebesgue integrable as well as the operator-valued function $\psi(t-\eta,\eta)$ is continuous
in uniform topology about the variables $t$ and $\eta$ for $0\leq t\leq a$ and $0\leq\eta\leq t-\epsilon$, we know that
\begin{eqnarray*}
J_6&\leq&
\sup\limits_{\eta\in
[0,t'-\epsilon]}\|\psi(t''-\eta,\eta)-\psi(t'-\eta,\eta)\| C\int_{0}^{t'-\epsilon}\int_{0}^{\eta}(\eta-s)^{\gamma-1}
\psi_R(s)dsd\eta\\[6pt]
& &+C^2\int_{t'-\epsilon}^{t'}\int_{0}^{\eta}
[(t''-\eta)^{\alpha-1}+(t'-\eta)^{\alpha-1}](\eta-s)^{\gamma-1}\psi_R(s)dsd\eta\\[6pt]
&\leq&
\Big(\frac{1-\beta}{\gamma-\beta}\Big)^{1-\beta}\frac{C(t')^\gamma
\|\psi_R\|_{L^{\frac{1}{\beta}}[0,a]}}{\gamma}
\sup\limits_{\eta\in
[0,t'-\epsilon]}\|\psi(t''-\eta,\eta)-\psi(t'-\eta,\eta)\|\\[6pt]
& &+C^2\Gamma(\gamma)\int_{t'-\epsilon}^{t'}
[(t''-\eta)^{\alpha-1}I_\eta^\gamma \psi_R(\eta)+(t'-\eta)^{\alpha-1}I_\eta^\gamma \psi_R(\eta)]d\eta
\end{eqnarray*}

\begin{eqnarray*}
\rightarrow0\quad \textrm{as}\quad t''-t'\rightarrow 0\quad \textrm{and}\quad \epsilon\rightarrow 0.\qquad\qquad\qquad\qquad\qquad\quad
\end{eqnarray*}\vskip2mm
\noindent
As a result, $\|(\mathcal {Q}u)(t'')-(\mathcal {Q}u)(t')\|\rightarrow 0$
independently of $u\in \Omega_{R}$ as $t''-t'\rightarrow 0$, which
means that the operator $\mathcal {Q}:\Omega_{R}\to \Omega_{R}$ is equicontinuous.\vskip3mm

Next, we prove that $\mathcal {Q}:F\rightarrow F$ is a convex-power condensing operator, where $F=\overline{\textrm{co}}\,\mathcal {Q}(\Omega_{R})$ and
$\overline{\textrm{co}}$ means the closure of convex hull. Then one can
easily to verify that the operator $\mathcal {Q}$ maps $F$ into itself and $F\subset
C(I,E)$ is equicontinuous. Let
$u_0\in F$. In the following, we will prove that there exists a positive integer
$n_0$ such that for any bounded and nonprecompact subset $D\subset
F$
\begin{equation} \label{39}
\mu_C\Big( \mathcal {Q}^{(n_0,u_0)}(D)\Big)<\mu_C(D).
\end{equation}
For any $D\subset F$ and $u_0\in F$, by the definition of operator $\mathcal {Q}^{(n,u_0)}$ and the equicontinuity
of $F$, we get that $\mathcal {Q}^{(n,u_0)}(D)\subset \Omega_R$ is also
equicontinuous. Therefore, we know from Lemma \ref{Le5} that
\begin{equation} \label{310}
\mu_C\Big(\mathcal {Q}^{(n,u_0)}(D)\Big)=\max\limits_{t\in I}\mu\Big(\mathcal {Q}^{(n,u_0)}
(D)(t)\Big),\qquad n=1,2,\cdots.
\end{equation}
 By Lemma \ref{Le3}, there
exists a countable set $D_1=\{u^1_n\}\subset D$, such that
\begin{equation} \label{311}
\mu(\mathcal {Q}(D)(t))\leq 2\mu(\mathcal {Q}(D_1)(t)).
\end{equation}
By the fact
$$
\int_0^au(s)ds\in a\overline{\textrm{co}}\{u(s)\mid s\in I\},\qquad \forall\;u\in C(I,E),
$$
we know that
\begin{equation} \label{312}
\mu\Big(\Big\{\int_0^tK(t,s)u(s)ds\mid u\in D,\;t
\in I\Big\}\Big)\leq aK_0\mu(\{u(t)\mid u\in D,\;t\in
I\})
\end{equation}
 and
\begin{equation} \label{313}
\mu\Big(\Big\{\int_0^aH(t,s)u(s)ds\mid u\in D,\;t\in I\Big\}
\Big)\leq aH_0\mu(\{u(t)\mid u\in D,\;t\in I\}).
\end{equation}
where $K_0=\max\limits_{(t,s)\in \triangle}|K(t,s)|$,
$H_0=\max\limits_{(t,s)\in \triangle_0}|H(t,s)|$. Therefore, by
\eqref{31}, \eqref{311}, \eqref{312}, \eqref{313}, Lemma \ref{Le1}, Lemma \ref{Le22}, Lemma \ref{Le4} and the assumption (F2), we have
\begin{eqnarray*}
& &\mu(\mathcal {Q}^{(1,u_0)}(D)(t))\\&= &\mu(\mathcal {Q}(D)(t))\leq
2\mu(\mathcal {Q}(D_1)(t))\qquad\qquad\qquad\qquad\qquad\qquad\qquad\qquad\quad
\end{eqnarray*}

\begin{equation} \label{314}\begin{split}
&\leq
2\mu\Big(A^{-1}(0)u_0+\int_0^t\psi(t-\eta,\eta)U(\eta)u_0d\eta\Big)\\[6pt]
&\quad+ 2\mu\Big(\Big\{\int_0^t\psi(t-\eta,\eta)
f(\eta,u_n^1(\eta),(Tu_n^1)(\eta),(Su_n^1)(\eta))d\eta\Big\}\Big)\\[6pt]
&\quad+ 2\mu\Big(\Big\{\int_0^t\int_0^\eta\psi(t-\eta,\eta)\varphi(\eta,s)\\[6pt]
&\quad\cdot
f(s,u_n^1(s),(Tu_n^1)(s),(Su_n^1)(s))dsd\eta\Big\}\Big)\\[6pt]&\leq 4C\int_{0}^{t}
(t-\eta)^{\alpha-1}[L_1\mu(D_1(\eta))+L_2\mu((TD_1)(\eta))\\[6pt]
&\quad+L_3\mu((SD_1)(\eta))]d\eta\\[6pt]
&\quad+ 8C^2\int_0^t\int_0^\eta(t-\eta)^{\alpha-1}(\eta-s)^{\gamma-1}\\[6pt]
&\qquad\cdot [L_1\mu(D_1(s))+L_2\mu((TD_1)(s))+L_3\mu((SD_1)(s))]dsd\eta\\[6pt]&\leq 4C\int_{0}^{t}
(t-\eta)^{\alpha-1}(L_1+aK_0L_2+aH_0L_3)\mu(D_1(\eta))d\eta\\[6pt]
&\quad+ 8C^2\mathbf{B}(\alpha,\gamma)\int_0^t(t-\eta)^{\alpha+\gamma-1}(L_1+aK_0L_2+aH_0L_3)
\mu(D_1(\eta))d\eta\\[6pt]
&\leq
\frac{4CM\Gamma(\alpha)t^\alpha}{\Gamma(1+\alpha)}\mu_C(D)
+\frac{8C^2M\Gamma(\alpha)\Gamma(\gamma)t^{\alpha+\gamma}}{\Gamma(1+\gamma+\alpha)}\mu_C(D),
\end{split}
\end{equation}
where
\begin{equation} \label{315}
M:=L_1+aK_0L_2+aH_0L_3.
\end{equation}

Again by Lemma \ref{Le3}, there exists a countable set
$D_2=\{u_n^2\}\subset \overline{\textrm{co}}\{\mathcal {Q}^{(1,u_0)}$ $ (D),
u_0\}$, such that
\begin{equation} \label{316}
\mu\Big(\mathcal {Q}(\overline{\textrm{co}}\{\mathcal {Q}^{(1,u_0)}(D),
u_0\})(t)\Big)\leq 2\mu(\mathcal {Q}(D_2)(t)).
\end{equation} Therefore, by
\eqref{31}, \eqref{312}, \eqref{313}, \eqref{314}, \eqref{315}, \eqref{316}, Lemma \ref{Le1}, Lemma \ref{Le22},  Lemma \ref{Le2} (iii), Lemma \ref{Le4} and the assumption (F2), we get that
\begin{equation} \label{317}\begin{split}
& \mu(\mathcal {Q}^{(2,u_0)}(D)(t))\\&= \mu\Big(\mathcal {Q}(\overline{\textrm{co}}\{\mathcal {Q}^{(1,u_0)}(D),
u_0\})(t)\Big)\\[6pt]&\leq
2\mu(\mathcal {Q}(D_2)(t))\\[6pt]&\leq
2\mu\Big(A^{-1}(0)u_0+\int_0^t\psi(t-\eta,\eta)U(\eta)u_0d\eta\Big)\\[6pt]&\quad+ 2\mu\Big(\Big\{\int_0^t\psi(t-\eta,\eta)
f(\eta,u_n^2(\eta),(Tu_n^2)(\eta),(Su_n^2)(\eta))d\eta\Big\}\Big)\\[6pt]
&\quad+ 2\mu\Big(\Big\{\int_0^t\int_0^\eta\psi(t-\eta,\eta)\varphi(\eta,s)
f(s,u_n^2(s),(Tu_n^2)(s),(Su_n^2)(s))dsd\eta\Big\}\Big)\\[6pt]&\leq 4CM\int_{0}^{t}
(t-\eta)^{\alpha-1}\mu(D_2(\eta))d\eta+ 8C^2M\mathbf{B}(\alpha,\gamma)\\[6pt]&\qquad \cdot\int_0^t(t-\eta)^{\alpha+\gamma-1}
\mu(D_2(\eta))d\eta\\[6pt]&
\leq 4CM\int_{0}^{t}
(t-\eta)^{\alpha-1}\mu\Big(\overline{\textrm{co}}\{\mathcal {Q}^{(1,u_0)}(D),
u_0\}(\eta)\Big)d\eta\\[6pt]&\quad+ 8C^2M\mathbf{B}(\alpha,\gamma)\int_0^t(t-\eta)^{\alpha+\gamma-1}
\mu\Big(\overline{\textrm{co}}\{\mathcal {Q}^{(1,u_0)}(D),
u_0\}(\eta)\Big)d\eta\\[6pt]&
\leq 4CM\int_{0}^{t}
(t-\eta)^{\alpha-1}\Big[\frac{4CM\Gamma(\alpha)t^\alpha}{\Gamma(1+\alpha)}
+\frac{8C^2M\Gamma(\alpha)\Gamma(\gamma)t^{\alpha+\gamma}}{\Gamma(1+\gamma+\alpha)}\Big]
d\eta\cdot\mu_C(D)\\[6pt]&\quad+ 8C^2M\mathbf{B}(\alpha,\gamma)\int_0^t(t-\eta)^{\alpha+\gamma-1}
\Big[\frac{4CM\Gamma(\alpha)t^\alpha}{\Gamma(1+\alpha)}\\[6pt]&\quad
+\frac{8C^2M\Gamma(\alpha)\Gamma(\gamma)t^{\alpha+\gamma}}
{\Gamma(1+\gamma+\alpha)}\Big]d\eta\cdot\mu_C(D)\\[6pt]
&=
\frac{4CM\cdot4CM\Gamma^2(\alpha)t^{2\alpha}}{\Gamma(1+2\alpha)}\mu_C(D)
+\frac{2\cdot4CM\cdot8C^2M\Gamma^2(\alpha)\Gamma(\gamma)t^{2\alpha+\gamma}}
{\Gamma(1+\gamma+2\alpha)}\mu_C(D)
\\[6pt]&\quad+\frac{8C^2M\cdot8C^2M\Gamma^2(\alpha)
\Gamma^2(\gamma)t^{2\alpha+2\gamma}}{\Gamma(1+2\gamma+2\alpha)}\mu_C(D).
\end{split}
\end{equation}

If for $\forall\;t\in I$, we assume that
\begin{equation} \label{318}\begin{split}
& \mu\Big(\mathcal {Q}^{(k,u_0)}(D)(t)\Big)\\[6pt]&\leq
\frac{C_k^0(4CM\Gamma(\alpha)t^{\alpha})^k\cdot(8C^2M\Gamma(\alpha)
\Gamma(\gamma)t^{\alpha+\gamma})^0}{\Gamma(1+k\alpha)}\mu_C(D)\\[6pt]&\quad
+\frac{C_k^1(4CM\Gamma(\alpha)t^{\alpha})^{k-1}\cdot(8C^2M\Gamma(\alpha)
\Gamma(\gamma)t^{\alpha+\gamma})^1}{\Gamma(1+\gamma+k\alpha)}\mu_C(D)\\[6pt]&\quad\;\,\vdots
\\[6pt]&\quad+\frac{C_k^{k-1}(4CM\Gamma(\alpha)t^{\alpha})^{1}\cdot(8C^2M\Gamma(\alpha)
\Gamma(\gamma)t^{\alpha+\gamma})^{k-1}}{\Gamma(1+(k-1)\gamma+k\alpha)}\mu_C(D)
\\[6pt]&\quad+\frac{C_k^{k}(4CM\Gamma(\alpha)t^{\alpha})^{0}\cdot(8C^2M\Gamma(\alpha)
\Gamma(\gamma)t^{\alpha+\gamma})^{k}}{\Gamma(1+k\gamma+k\alpha)}\mu_C(D).
\end{split}\end{equation}
Then by Lemma \ref{Le3}, there exists a countable set
$D_{k+1}=\{u_n^{k+1}\}\subset
\overline{\textrm{co}}\{\mathcal {Q}^{(k,u_0)}$ $(D), u_0\}$, such that
\begin{equation} \label{319}
\mu\Big(\mathcal {Q}(\overline{\textrm{co}}\{\mathcal {Q}^{(k,u_0)}(D),
u_0\})(t)\Big)\leq 2\mu(\mathcal {Q}(D_{k+1})(t)).
\end{equation} From \eqref{31}, \eqref{312}, \eqref{313}, \eqref{315}, \eqref{318}, \eqref{319}, Lemma \ref{Le1}, Lemma \ref{Le22}, Lemma \ref{Le2} (iii), Lemma \ref{Le4}, and the assumption (F2) and proper integral transformation, we get that
\begin{eqnarray*}
& &  \mu(\mathcal {Q}^{(k+1,u_0)}(D)(t))\\[6pt]&=&\mu\Big(\mathcal {Q}(\overline{\textrm{co}}\{\mathcal {Q}^{(k,u_0)}(D),
u_0\})(t)\Big)\leq
2\mu(\mathcal {Q}(D_{k+1})(t))\\[6pt]&\leq&
4CM\int_{0}^{t}
(t-\eta)^{\alpha-1}\mu(D_{k+1}(\eta))d\eta\\[6pt]&\quad& + 8C^2M\mathbf{B}(\alpha,\gamma)\int_0^t(t-\eta)^{\alpha+\gamma-1}
\mu(D_{k+1}(\eta))d\eta\\[6pt]&
\leq&  4CM\int_{0}^{t}
(t-\eta)^{\alpha-1}\mu\Big(\overline{\textrm{co}}\{\mathcal {Q}^{(k,u_0)}(D),
u_0\}(\eta)\Big)d\eta\\[6pt]&\quad& + 8C^2M\mathbf{B}(\alpha,\gamma)\int_0^t(t-\eta)^{\alpha+\gamma-1}
\mu\Big(\overline{\textrm{co}}\{\mathcal {Q}^{(k,u_0)}(D),
u_0\}(\eta)\Big)d\eta
\end{eqnarray*}
\begin{equation} \label{320}\begin{split}
&\leq  \frac{C_{k+1}^0(4CM\Gamma(\alpha)t^{\alpha})^{k+1}\cdot(8C^2M\Gamma(\alpha)
\Gamma(\gamma)t^{\alpha+\gamma})^0}{\Gamma(1+(k+1)\alpha)}\mu_C(D)\\[8pt]&\quad
+\frac{C_{k+1}^1(4CM\Gamma(\alpha)t^{\alpha})^{k}\cdot(8C^2M\Gamma(\alpha)
\Gamma(\gamma)t^{\alpha+\gamma})^1}{\Gamma(1+\gamma+(k+1)\alpha)}\mu_C(D)\\[8pt]&\quad \;\,\vdots
\\[8pt]&\quad +\frac{C_{k+1}^{k}(4CM\Gamma(\alpha)t^{\alpha})^{1}\cdot(8C^2M\Gamma(\alpha)
\Gamma(\gamma)t^{\alpha+\gamma})^{k}}{\Gamma(1+k\gamma+(k+1)\alpha)}\mu_C(D)
\\[8pt]&\quad +\frac{C_{k+1}^{k+1}(4CM\Gamma(\alpha)t^{\alpha})^{0}\cdot(8C^2M\Gamma(\alpha)
\Gamma(\gamma)t^{\alpha+\gamma})^{k+1}}{\Gamma(1+(k+1)\gamma+(k+1)\alpha)}\mu_C(D).
\end{split}\end{equation}
Therefore, by the method of mathematical induction, we know that for any positive
integer $n$ and $t\in I$
\begin{equation} \label{321}\begin{split}
& \mu\Big(\mathcal {Q}^{(n,u_0)}(D)(t)\Big)\\[8pt]&\leq
\frac{C_n^0(4CM\Gamma(\alpha)t^{\alpha})^n\cdot(8C^2M\Gamma(\alpha)
\Gamma(\gamma)t^{\alpha+\gamma})^0}{\Gamma(1+n\alpha)}\mu_C(D)\\[8pt]&\quad
+\frac{C_n^1(4CM\Gamma(\alpha)t^{\alpha})^{n-1}\cdot(8C^2M\Gamma(\alpha)
\Gamma(\gamma)t^{\alpha+\gamma})^1}{\Gamma(1+\gamma+n\alpha)}\mu_C(D)\\[8pt]&\quad\;\,\vdots
\\[8pt]&\quad+\frac{C_n^{n-1}(4CM\Gamma(\alpha)t^{\alpha})^{1}\cdot(8C^2M\Gamma(\alpha)
\Gamma(\gamma)t^{\alpha+\gamma})^{n-1}}{\Gamma(1+(n-1)\gamma+n\alpha)}\mu_C(D)
\\[8pt]&\quad+\frac{C_n^{n}(4CM\Gamma(\alpha)t^{\alpha})^{0}\cdot(8C^2M\Gamma(\alpha)
\Gamma(\gamma)t^{\alpha+\gamma})^{n}}{\Gamma(1+n\gamma+n\alpha)}\mu_C(D).
\end{split}\end{equation}
Hence, by \eqref{310} and
\eqref{321}, we get that

\begin{eqnarray*}
\mu_C\Big(\mathcal {Q}^{(n,u_0)}(D)\Big)&=&\max\limits_{t\in
I}\mu\Big(\mathcal {Q}^{(n,u_0)}(D)(t)\Big)\\[6pt]&\leq&
\Big[\frac{(4CM\Gamma(\alpha))^n\cdot a^{n\alpha}}{\Gamma(1+n\alpha)}\\[8pt]&\quad&
+\frac{C_n^1(4CM\Gamma(\alpha))^{n-1}\cdot(8C^2M\Gamma(\alpha)
\Gamma(\gamma))^1\cdot a^{\gamma+n\alpha}}{\Gamma(1+\gamma+n\alpha)}\qquad\qquad\qquad\qquad\quad
\end{eqnarray*}

\begin{equation} \label{322}\begin{split}
\qquad\qquad\quad&\quad\;\,\vdots
\\[8pt]&\quad+\frac{C_n^{n-1}(4CM\Gamma(\alpha))^{1}\cdot(8C^2M\Gamma(\alpha)
\Gamma(\gamma))^{n-1}\cdot a^{(n-1)\gamma+n\alpha}}{\Gamma(1+(n-1)\gamma+n\alpha)}
\\[8pt]&\quad+\frac{(8C^2M\Gamma(\alpha)
\Gamma(\gamma))^{n}\cdot a^{n\gamma+n\alpha}}{\Gamma(1+n\gamma+n\alpha)}\Big]\mu_C(D).
\end{split}\end{equation}
Thanks to the well-known Stirling's formula
$$\Gamma(1+s)=\sqrt{2\pi s}\Big(\frac{s}{e}\Big)^{s}e^{\frac{\theta}{12s}},\quad s>0,\;0<\theta <1,$$
we get that when $n\rightarrow \infty$, then
\begin{eqnarray*}
\frac{(4CM\Gamma(\alpha))^n\cdot a^{n\alpha}}{\Gamma(1+n\alpha)}=\frac{(4CM\Gamma(\alpha))^n\cdot a^{n\alpha}} {\sqrt{2\pi
n\alpha}(\frac{n\alpha}{e})^{n\alpha}e^{\frac{\theta}{12n\alpha}}}
&\rightarrow& 0,
\end{eqnarray*}\vskip2mm
\begin{eqnarray*}
& &\frac{C_n^1(4CM\Gamma(\alpha))^{n-1}\cdot(8C^2M\Gamma(\alpha)
\Gamma(\gamma))^1\cdot a^{\gamma+n\alpha}}{\Gamma(1+\gamma+n\alpha)}\\[8pt] &\qquad&=\frac{C_n^1(4CM\Gamma(\alpha))^{n-1}\cdot(8C^2M\Gamma(\alpha)
\Gamma(\gamma))^1\cdot a^{\gamma+n\alpha}} {\sqrt{2\pi
(\gamma+n\alpha)}(\frac{\gamma+n\alpha}{e})^{\gamma+n\alpha}
e^{\frac{\theta}{12(\gamma+n\alpha)}}}\\[10pt]
&\qquad&\rightarrow 0,\qquad\qquad\qquad\qquad\qquad\qquad\qquad\qquad\qquad\qquad
\end{eqnarray*}
$$
\cdots \qquad    \cdots\qquad\qquad\qquad\qquad\qquad\qquad\quad
$$
\begin{eqnarray*}
\\& &\frac{C_n^{n-1}(4CM\Gamma(\alpha))^{1}\cdot(8C^2M\Gamma(\alpha)
\Gamma(\gamma))^{n-1}\cdot a^{(n-1)\gamma+n\alpha}}{\Gamma(1+(n-1)\gamma+n\alpha)}\\[8pt] &\qquad&=\frac{C_n^{n-1}(4CM\Gamma(\alpha))^{1}\cdot(8C^2M\Gamma(\alpha)
\Gamma(\gamma))^{n-1}\cdot a^{(n-1)\gamma+n\alpha}} {\sqrt{2\pi
((n-1)\gamma+n\alpha)}(\frac{(n-1)\gamma+n\alpha}{e})^{(n-1)\gamma+n\alpha}
e^{\frac{\theta}{12((n-1)\gamma+n\alpha)}}}\\[10pt]
&\qquad&\rightarrow 0
\end{eqnarray*}
and\vskip2mm
\begin{eqnarray*}
\frac{(8C^2M\Gamma(\alpha)
\Gamma(\gamma))^{n}\cdot a^{n\gamma+n\alpha}}{\Gamma(1+n\gamma+n\alpha)}&=&\frac{(8C^2M\Gamma(\alpha)
\Gamma(\gamma))^{n}\cdot a^{n\gamma+n\alpha}} {\sqrt{2\pi
(n\gamma+n\alpha)}(\frac{n\gamma+n\alpha}{e})^{n\gamma+n\alpha}
e^{\frac{\theta}{12(n\gamma+n\alpha)}}}\\[6pt]
& \rightarrow& 0.\qquad\qquad\qquad\qquad\qquad\qquad\qquad\qquad\qquad\qquad\qquad\,
\end{eqnarray*}

Therefore, there must exist a positive integer $n_0$, which is large enough,  such
that
\begin{equation} \label{323}
\frac{(4CM\Gamma(\alpha))^{n_0}\cdot a^{n_0\alpha}}{\Gamma(1+n_0\alpha)}
+\cdots+\frac{(8C^2M\Gamma(\alpha)
\Gamma(\gamma))^{n_0}\cdot a^{n_0\gamma+n_0\alpha}}{\Gamma(1+n_0\gamma+n_0\alpha)}<1.
\end{equation}
Hence, from \eqref{322} and
\eqref{323} we know that \eqref{39} is satisfied, which means that $\mathcal {Q}:F\rightarrow F$ is a convex-power condensing operator. It
follows from Lemma \ref{Le6} that the operator $\mathcal {Q}$ defined by \eqref{31} has at least one fixed point $u\in
F$, which is just a mild solution of initial value problem \eqref{11}. This completes
the proof of Theorem
\ref{th1}. \hfill$\Box$

\vskip3mm
\noindent\textbf{Proof of Theorem \ref{th2}.}\quad From the proof of Theorem \ref{th1}, we know that the mild solution of initial value problem  \eqref{11} is equivalent to the fixed point of the operator $\mathcal {Q}$ defined by \eqref{31}. In what follows, we prove that there exists a positive constant $R$ such that the operator $\mathcal {Q}$
maps the set $\Omega_R$ to $\Omega_R$. For any $u\in \Omega_R$ and a.e. $t\in I$, by \eqref{31}, \eqref{2}, Lemmas
\ref{Le1}-\ref{Le22}, the assumption (F1)$^\ast$ and H\"{o}lder inequality, we know that
\begin{equation} \label{324}\begin{split}
\|(\mathcal {Q}u)(t)\|&\leq \|A^{-1}(0)u_0\|+\Big\|\int_0^{t}
\psi({t}-\eta,\eta)U(\eta)g(u)d\eta\Big\|\\[4pt]
&\quad +\Big\|\int_0^{t }\psi({t }-\eta,\eta)
f(\eta,u (\eta),(Tu )(\eta),(Su_u)(\eta))d\eta\Big\|\\[4pt]
&\quad +\Big\|\int_0^{t }\int_0^\eta\psi({t }-\eta,\eta)\varphi(\eta,s)
f(s,u (s),(Tu )(s),(Su )(s))dsd\eta\Big\|\\[4pt]
&\leq C\|u_0\|+C^2\int_0^{t }({t }-\eta)^{\alpha-1}(1+\eta^\gamma)\|u_0\|d\eta\\[4pt]
&\quad +C\int_0^{t }({t }-\eta)^{\alpha-1}\phi (\eta)\Phi(\|u\|)d\eta\\[4pt]
&\quad
+C^2\int_0^{t }\int_0^\eta ({t }-\eta)^{\alpha-1}(\eta-s)^{\gamma-1}\phi (s)\Phi(\|u\|)dsd\eta\\[4pt]
&\leq C\|u_0\|+C^2\|u_0\|{t }^\alpha\Big(\frac{1}{\alpha}+({t })^\gamma \mathbf{B}(\alpha,\gamma+1)\Big)\\[4pt]
&\quad +C\Phi(R)\int_0^{t }({t }-\eta)^{\alpha-1}\phi (\eta)d\eta\\[4pt]&\quad
+C^2\Phi(R)\mathbf{B}(\alpha,\gamma)\int_0^{t } ({t }-\eta)^{\alpha+\gamma-1}\phi (\eta)d\eta\qquad\qquad\qquad\qquad\qquad\qquad\qquad\qquad\qquad\qquad\qquad\qquad
\end{split}
\end{equation}
\begin{eqnarray*}
&\leq& C\|u_0\|\delta_2+C\Phi(R)\Big(\int_0^{t }({t }-\eta)^{\frac{\alpha-1}{1-\beta}}d\eta\Big)^{1-\alpha_1} \Big(\int_0^{t }\phi ^{\frac{1}{\beta}}(\eta)d\eta\Big)^{\beta}\\[4pt]
&\quad&+C^2\Phi(R)\mathbf{B}(\alpha,\gamma)\Big(\int_0^{t }({t }-\eta)^{\frac{\alpha+\gamma-1}
{1-\beta}}d\eta\Big)^{1-\beta} \Big(\int_0^{t }\phi ^{\frac{1}{\beta}}(\eta)d\eta\Big)^{\beta}\\[4pt]
&\leq& C\|u_0\|\delta_2+ C\delta_1\Phi(R)a^{\alpha-\beta}\|
\phi \|_{L^{\frac{1}{\beta}}[0,a]}\\[4pt]
&\leq& R.\qquad\qquad\qquad\qquad\qquad\qquad\qquad\qquad\;
\end{eqnarray*}
Therefore, from \eqref{324} we get that the operator $\mathcal {Q}:\Omega_{R}\to \Omega_{R}$. By adopting a completely similar method with which used in the proof
of Theorem \ref{th1}, we can prove that $\mathcal {Q}:\Omega_{R}\to \Omega_{R}$ is continuous and equicontinuous, and also $\mathcal {Q}:\overline{\textrm{co}}\,\mathcal {Q}(\Omega_{R})\to \overline{\textrm{co}}\,\mathcal {Q}(\Omega_{R})$ is a convex-power condensing operator. By Lemma \ref{Le6}, we know that the operator $\mathcal {Q}$ defined by \eqref{31} has at least one fixed point $u\in
\overline{\textrm{co}}\,\mathcal {Q}(\Omega_{R})$, which is just a mild solution of initial value problem \eqref{11}. This completes
the proof of Theorem
\ref{th2}. \hfill$\Box$
\vskip3mm

\begin{remark}\label{re31}
From \eqref{310} and \eqref{314} of Theorem \ref{th1} we know that
if we assume
\begin{equation} \label{325}
\frac{4CM\Gamma(\alpha)a^\alpha}{\Gamma(1+\alpha)}
+\frac{8C^2M\Gamma(\alpha)\Gamma(\gamma)a^{\alpha+\gamma}}{\Gamma(1+\gamma+\alpha)}<1
\end{equation}
directly, we can
apply the famous Sadoveskii's fixed point theorem to obtain the results  of Theorem \ref{th1} and Theorem \ref{th2}. However, from the above arguments, one can find that
we do not need this redundant condition \eqref{325} by virtue of fixed point theorem with respect to convex-power condensing operator.
\end{remark}

\vskip 10mm
\section{Application}
In this section, we present an example, which do not aim at generality but indicate how our
abstract result can be applied to concrete problem. As an application, we consider the following initial value problem of time fractional non-autonomous partial differential equation with homogeneous Dirichlet boundary condition
\begin{equation}\label{41}
\small\left\{\begin{array}{ll}
   \frac{\partial^\alpha}{\partial t^\alpha}u(x,t)-\kappa(x,t)\Delta u(x,t)
=\frac{\sin(\pi t)}{1+| u(x,\,t)|}+e^{-t}\sin\Big(\int_0^t(t-s)u(x,s)ds\Big)
\\[6pt] \qquad \qquad \qquad \qquad \qquad \qquad \quad +e^{-t}\cos\Big(\int_0^1e^{-|t-s|}u(x,s)ds\Big),\; x\in\Omega,\; t\in I,
\\[8pt]
u(x,t)=0,\qquad x\in \partial\Omega,\quad t\in I,\\[8pt]
u(x,0)=(\kappa(x,0))^{-1}\varphi(x),\quad x\in\Omega,
\end{array} \right.
\end{equation}\vskip 0.1mm
\noindent where  $\frac{\partial^\alpha}{\partial t^\alpha}$ is the Caputo fractional
order partial derivative of order $\alpha$, $0<\alpha\leq1$, $I=[0,1]$, $\Delta$ is the Laplace operator, $\Omega\subseteq \mathbb{R}^n$ is a bounded domain with a sufficiently smooth boundary $\partial\Omega$, the coefficient of heat conductivity $\kappa(x,t)$ is continuous on $\Omega\times[0,1]$ and it is is uniformly H\"{o}lder continuous
in $t$, which means that for any $t_1$, $t_2\in I$, there exist a constant $0<\gamma\leq1$ and a positive constant $C$ independent of $t_1$ and $t_2$, such that
\begin{equation} \label{42}
|\kappa(x,t_2)-\kappa(x,t_1)|\leq C|t_2-t_1|^\gamma,\quad x\in \Omega,
\end{equation}
and $\varphi\in L^2(\Omega)$.\vskip 2mm

Let $E=L^2(\Omega)$ be a Banach space with the $L^2$-norm $\|\cdot\|_2$. We define an operator $A(t)$ in Banach space $E$ by
\begin{equation} \label{44}
D(A)=H^2(\Omega)\cap H_0^1(\Omega),\qquad A(t)u=-\kappa(x,t)\Delta u,
\end{equation}
where $H^2(\Omega)$ is the completion of the space $C^2(\Omega)$ with respect to the norm
$$
\|u\|_{H^2(\Omega)}=\Big(\int_\Omega\sum_{|\vartheta|\leq 2}|D^\vartheta u(x)|^2dx\Big)^{\frac{1}{2}},
$$
$C^2(\Omega)$ is the set of all continuous functions defined on $\Omega$ which have continuous partial derivatives of order less than or equal to $2$, $H_0^1(\Omega)$ is the completion of $C^1(\Omega)$ with
respect to the norm $\|u\|_{H^1(\Omega)}$, and $C_0^1(\Omega)$ is the set of all functions $u\in C^1(\Omega)$ with compact
supports on the domain $\Omega$. Then it is well known from \cite{fr69} that $-A(s)$ generates an analytic semigroup $e^{-tA(s)}$ in $E$. By \eqref{42} and \eqref{44} one can easily to verify that the linear operator $-A(t)$ satisfies the assumptions (A1) and (A2). \vskip 3mm

Further, for any $t\in [0,1]$, we define
$$u(t)=u(\cdot,t),\quad K(t,s)=t-s\quad \textrm{for} \quad0\leq s\leq t\leq1,$$
$$
a=1,\quad H(t,s)=e^{-|t-s|}\quad \textrm{for} \quad0\leq s,\, t\leq1,
$$
$$
(Tu)(t)=\int_0^tK(t,s)u(\cdot,s)ds,
\quad (Su)(t)=\int_0^1H(t,s)u(\cdot,s)ds,
$$
$$
f(t,u(t),(Tu)(t),(Su)(t))=\frac{\sin(\pi t)}{1+| u(\cdot,\,t)|}+e^{-t}\sin\big((Tu)(t)\big)+e^{-t}\cos\big((Su)(t)\big),
$$
$$
A^{-1}(0)=(\kappa(\cdot,0))^{-1},\quad u_0=\varphi(\cdot).
$$
Then the initial value problem of
time fractional non-autonomous partial differential equation with homogeneous Dirichlet boundary condition \eqref{41} can be transformed into the abstract form of initial value problem to time fractional non-autonomous
integro-differential evolution equation of mixed type \eqref{11}.\vskip 3mm

\begin{theorem}\label{th41}
The initial value problem of time fractional non-autonomous partial differential equation with homogeneous Dirichlet boundary condition \eqref{41} has at least one mild solution $u\in C(\Omega\times [0,1])$.
\end{theorem}\vskip 2mm
\proof By the definition of nonlinear term $f$ one can easily to verify that the assumption (F1) is satisfied with
$$
\psi_r(t)=\sqrt{\textrm{mes}\Omega}(\sin(\pi t)+2e^{-t}),\quad \textrm{and} \quad \beta=\rho=0.
$$
From the definition of nonlinear term $f$, we know that $f(t,u,v,w)$ is Lipschitz continuous about the variables $u$, $v$ and $w$ with Lipschitz constants $k_u=1$, $k_v=1$ and $k_w=1$, respectively. Therefore, by Lemma \ref{Le2} (8) we get that the assumption (F2) is satisfied with positive constants
$$
L_1=L_2=L_3=1.
$$
From the fact  $\rho=0$ one can easily to verify that the condition \eqref{1} holds. Therefore, all the assumptions of Theorem \ref{th1} are satisfied. Hence,
the initial value problem of time fractional non-autonomous partial differential equation with homogeneous Dirichlet boundary condition \eqref{41} has at least one mild solution $u\in C(\Omega\times [0,1])$ duo to Theorem \ref{th1}.  This completes the proof of Theorem \ref{th41}.\endproof
\vskip0.1mm
\subsection*{Acknowledgements}
This work is supported by National Natural Science Foundations of China (No. 11501455, No. 11661071) and Doctoral Research Fund of Northwest Normal University (No. 6014/0002020209).\vskip15mm

\bibliographystyle{abbrv}

\vskip3mm

\end{document}